\documentclass[11pt]{tran-l}
\usepackage{amssymb}
 \newcommand{\Z}{{\mathbb Z}}
 \newcommand{\C}{{\mathbb C}}
\newcommand{\F}{{\mathcal F}}
 \newcommand{\Q}{{\mathbb Q}}
 \newcommand{\R}{{\mathbb R}}
 \newcommand{\K}{{\mathbb K}}

\renewcommand{\P}{{\mathbb P}}

\newcommand{\g}{\mathfrak g}

\renewcommand{\o}{\overline}
 \newtheorem{theorem}{Theorem}
 \newtheorem{lemma}[theorem]{Lemma} 
  \newtheorem{proposition}[theorem]{Proposition}
 \newtheorem{corollary}[theorem]{Corollary}

 \newtheorem{definition}[theorem]{Definition}   
 \newtheorem{conjecture}{Conjecture}

 \newtheorem{fact}[theorem]{Fact}
\newcommand{\mo}{{\mathcal{O}}}

 \def\Box
  {\hfill \thinspace\vbox{\hrule height .5pt \hbox{\vrule  
   width .5pt \vbox to 7pt{\hbox to 3.5pt{}} \vrule width .5pt} 
   \hrule height 0pt depth .5pt}}

\begin{document}
\title{Cut numbers of $3$-manifolds}

\author{Adam S. Sikora}
\address{Dept. of Mathematics, 244 Math. Bldg., SUNY at Buffalo, 
Buffalo, NY 14260}

\email{asikora@buffalo.edu}
\subjclass{57M05, 57M27, 20F34,11E76} 
\keywords{cut number, 3-manifold, corank, skew-symmetric form, 
cohomology ring}

\date{}

\begin{abstract} 
We investigate the relations between the cut number, $c(M),$ and
the first Betti number, $b_1(M),$ of $3$-manifolds $M.$
We prove that the cut number of a ``generic'' $3$-manifold
$M$ is at most $2.$ 
This is a rather unexpected result since 
specific examples of $3$-manifolds
with large $b_1(M)$ and $c(M)\leq 2$ are hard to construct.
We also prove that for any complex semisimple Lie algebra $\g$ there exists
a $3$-manifold $M$ with $b_1(M)=dim\, \g$ and $c(M)\leq rank\, \g.$ 
Such manifolds can be explicitly constructed.
\end{abstract}

\maketitle

\section{Introduction}

Let $M$ be a closed, smooth manifold. The cut number of $M,$ 
$c(M),$ is the largest number of disjoint two-sided hypersurfaces 
$\F_1,...,\F_c,$ which do not separate $M,$ ie. surfaces such that
$M\setminus (\F_1\cup ...\cup \F_c)$ is connected.

By \cite[Prop. 4.2]{St} (see also our Proposition \ref{cut=corank}), 
$c(M)$ is the corank of 
$\pi_1(M).$ The corank of a group $\Gamma$ is the largest
number $n$ for which the free group on $n$
generators, $F_n,$ is an epimorphic image of $\pi_1(M).$ 

Although the corank of a group is not a very tangible quantity,
it can be calculated (at least theoretically) for any finitely presented
group. This calculation reduces to finding the
ranks of solutions of certain systems of equations on free groups.
The algorithm for finding such ranks was given in \cite[\S 9]{Ra}.

Since $c(M)$ is bounded from above by the first Betti number of $M,$ $b_1(M),$
it is natural to ask whether there exists an analogous lower bound
on $c(M).$ In particular, T. Kerler and J. H. Przytycki asked whether
\begin{equation}\label{cnc}
\frac{1}{3}b_1(M)\leq c(M)
\end{equation}
for any closed, oriented $3$-manifold $M.$
This question was motivated
by the following two inequalities:
$$\frac{1}{3}b_1(M)\leq \mo_p(M)/\mo_p(S^2\times S^1),$$
and
$$c(M)\leq \mo_p(M)/\mo_p(S^2\times S^1),\ \text{for $p=5,$}$$
proved by Cochran and Melvin, \cite[Thm 4.3]{CM}, and Gilmer and Kerler
respectively. Here $\mo_p(M)$ denotes the quantum order of $M,$ that is
the largest $k$ such that
$\tau_p(M)/(\xi_p-1)^k\in \Z[\xi_p],$ where $\xi_p=e^{2\pi i/p}$
and $\tau_p(M)$ is the $SO(3)$-quantum invariant of 3-manifolds.
(We assume that $\tau_p(S^3)=1$). Note that (\ref{cnc}) 
holds for (a) Seifert $3$-manifolds, (b) $3$-manifolds with
abelian fundamental groups, (c) $3$-manifolds $M$ with
$b_1(M)\leq 3.$ Furthermore, for closed surfaces
$\F,$ $\frac{1}{2}b_1(\F)=c(\F).$

We are going to prove the following two results contradicting (\ref{cnc}):

\begin{theorem}\label{thm-lie}
For any semisimple Lie algebra $\g$ over $\C$
there exists a closed, oriented $3$-manifold $M$ with $b_1(M)=dim\, \g$ and 
$c(M)\leq rank\, \g.$
\end{theorem}

The proof and information on how to construct such $3$-manifolds
can be found in Section \ref{s_lie}.

The dimensions of simple Lie algebras of rank $r$ are
as follows:
$$dim\, sl_{r+1}=r(r+2),\quad dim\, so_{2r+1}=dim\, sp_{2r}=r(2r+1),\quad
dim\, so_{2r}=r(2r-1),$$ 
$$dim\, G_2=14,\ dim\, F_4=52,\ dim\, E_6=78,\ dim\, E_7=133,\ 
dim\, E_8=248.$$
Therefore, Theorem \ref{thm-lie} provides many counterexamples for 
(\ref{cnc}). In particular, it implies the existence of 
$3$-manifolds with
$b_1(M)=8,10, 14$ and $c(M)\leq 2.$

Let $\Omega_{\Q}(n)$ be the space of all skew-symmetric 
$3$-forms on $\Q^n,$ $\Omega_{\Q}(n)=(\bigwedge^3 \Q^n)^*.$ 
For a given closed, oriented $3$-manifold $M,$ a choice of a basis of 
$H^1(M,\Q)$ provides an isomorphism $H^1(M,\Q)=\Q^n$ for $n=b_1(M)$ and
makes it possible to consider the cup product $3$-form on $H^1(M,Q),$ 
$\Psi_M(x,y,z)=(x\cup y\cup z)[M],$ as an element of $\Omega_\Q(n).$

The cup $3$-form, $\Psi_M,$ determines the rational cohomology ring of 
$M.$ Since by a theorem of Sullivan, \cite{Su} (compare also \cite{T2}),
every $3$-form corresponds to the cup product form
for some closed oriented $3$-manifold, $\Omega_{\Q}(n)$ can be considered
as a space classifying the rational cohomology rings of $3$-manifolds $M$ with
a specified basis of $H^1(M,\Q).$

\begin{theorem}\label{main}
(1) The set of forms $\Psi_M\in \Omega_{\Q}(n)$ corresponding to 
$3$-manifolds $M$ with 
$c(M)\leq 2$ contains an open, dense subset 
(in the Zariski and Euclidean topologies) 
of $\Omega_{\Q}(n).$\\
(2) In particular, for any $n$ there exists a $3$-manifold $M$
with $b_1(M)=n$ and $c(M)\leq 2.$\\
(3) For any $n\geq 9$ there exist infinitely many $3$-manifolds $M$ with
$b_1(M)=n,$ $c(M)\leq 2,$ and with pairwise non-isomorphic rational 
cohomology rings.
\footnote{The assumption $n\geq 9$ is related to the fact that
for any $n\leq 8$ the cohomology rings, $H^*(M,\C),$
of $3$-manifolds $M$ with $b_1(M)=n$ belong to a finite number of isomorphism
classes only.}
\end{theorem}

Part (1) of the above result implies that the cut number
of a ``generic'' $3$-manifold is at most $2.$ This is a rather surprising
result, since it is not easy to construct specific examples of $3$-manifolds
with large $b_1(M)$ and $c(M)\leq 2.$

In an independent work on Kerler's and Przytycki's question,
S. Harvey found an explicit family of $3$-manifolds $M_n$ with $b_1(M_n)=n$
and $c(M_n)=1,$ \cite{Har}\footnote{Both S. Harvey and I announced
our results in September 2001.}. For the purpose of her 
construction, one takes
a link $L_n\subset S^3$ of $n$ components, obtained from the unlink by
linking every pair of components in the Whitehead link manner.
Harvey's manifold $M_n$ is the result of the $0$-surgery on $L_n.$

In yet another independent work inspired by the same problem, 
C. Leininger and A. Reid showed that there exist 
examples of a hyperbolic $3$-manifolds with $c(M)=1$ and $b_1(M)=5,$ \cite{LR}.

{\bf Acknowledgments:} We would like to thank S. Boyer, T. Kerler,
T. Goodwillie, J. Millson, J. H. Przytycki, J. Roberts, A. Razborov, 
G. Schwartz, J. Stallings, and V. Turaev for their helpful comments.

\section{Cut number and its siblings}
We will denote the corank of $\Gamma$ by $d(\Gamma).$

\begin{proposition}(cmp. \cite[Prop. 4.2]{St})\label{cut=corank}
For a closed, smooth manifold $M,$
$c(M)=d(\pi_1(M)).$
\end{proposition}

\begin{proof} Any epimorphism $\pi_1(M)\twoheadrightarrow F_c$
induces a map $f: M\twoheadrightarrow \bigvee_1^c S^1$
which can be assumed smooth outside the preimage of the base point.
If $y_1,...,y_c$ are regular values of $f$ lying in different
circles, then $f^{-1}(y_1),...,f^{-1}(y_c)$ are disjoint, two-sided 
surfaces in $M.$
Since these surfaces represent $c$ linearly independent elements of
$H_2(M,\Z/2),$ for each $i$ one can choose a connected component 
$\F_i\subset f^{-1}(y_i)$
such that $\F_1,...,\F_c$ are linearly independent in $H_2(M,\Z/2)$ as well.
The following lemma shows that $c(M)\geq d(\pi_1(M)).$

\begin{lemma} $\F_1,...,\F_c$ do not separate $M.$
\end{lemma} 

\begin{proof} Let loops $\gamma_1,...,\gamma_c$ represent elements
of $H_1(M,\Z/2)$ dual to $\F_1,...,\F_c,$
$\gamma_i\cdot \F_j=\delta_{ij}\ mod\ 2.$
Assume that $\gamma_i$ intersects $\F_j$ in at least two points $p_1,p_2.$ 
Let $C_1,C_2$ be arcs connecting $p_1$ and $p_2,$ such that
$C_1$ and $C_2$ lie on different sides of $\F_j$ and do not intersect other
surfaces.
By cutting $\gamma_i$ at $p_1$ and $p_2$ and inserting
the arcs $C_1$ and $C_2$ into $\gamma_i$ we obtain a new curve
with two less intersections with $\F_j$ and the number of 
intersections with $\F_{j'}$ unchanged for $j'\ne j.$
By applying this operation as many times as necessary,
we get closed curves $\gamma_1',...,\gamma_c'$ such that $\gamma_i'$
intersects $\F_i$ only once and it does not intersect $\F_j$ for 
$j\ne i.$ Consider a path $\gamma$ in $M$
joining two arbitrary points $x$ and $y$ in $M\setminus (\F_1,...,\F_c)$.
By by deforming $\gamma$ appropriately and 
pasting $\gamma\cdot \F_i$ copies of $\gamma_i'$ into $\gamma$
we obtain a new path connecting $x$ and $y$ which does not intersect 
$\F_1,...,\F_c.$
\end{proof}

We continue proving Proposition \ref{cut=corank}.
Let $\F_1,...,\F_c$ be disjoint, two-sided surfaces which do not
separate $M$ and let $N_1,.., N_c$ be disjoint neighborhoods
of $\F_1,...,\F_c$ provided with parameterizations 
$N_i\stackrel{\pi_i\times \lambda_i}{\longrightarrow} [-1,1]\times \F_i.$ 
Consider the map $f: M\to \bigvee_1^c S^1$ composed of the constant map
$M\setminus (N_1\cup...\cup N_c) \to *$ and the maps
$N_i\stackrel{\pi_i}\to [-1,1]\to S^1,$ where $S^1,$ obtained by identifying
the end points of $[-1,1],$
is the $i$-th circle in $\bigvee_1^c S^1.$ Since $f$
induces an epimorphism
$f_*: \pi_1(M) \twoheadrightarrow F_c,$ $d(\pi_1(M))\geq c(M).$
\end{proof}

The cut number can be also defined for simplicial
complexes:

\begin{definition}\label{simplical_cut}
For any simplicial complex $X$ define $c(X)$ to be the largest 
number $n$ such that there exists a subdivision $X'$ of $X$ and
$1$-cocycles $\alpha^1,...,\alpha^n\in C^1(X')$ such that
$\alpha^i\cup \alpha^j=0$ in $C^2(X')$ for all $i\ne j,$ and such that
$[\alpha^1],...,[\alpha^n] \in H^1(X,\Z)$
are linearly independent.
\end{definition}

Above, we assume that $X'$ is an ordered simplicial complex
and the cup product on $C^1(X')$ is given by the Alexander-Whitney 
diagonal approximation,
$$(\alpha\cup \beta)<x_0,x_1,x_2>=-\alpha(<x_0,x_1>)\beta(<x_1,x_2>),$$
compare \cite[VI\S 4]{Bre}.

The following generalizes Proposition 4.4 in \cite{Dw}.

\begin{proposition}\label{c=d} $c(X)=d(\pi_1(X))$
for any simplicial complex $X.$
\end{proposition}

\begin{proof}
If $d(\pi_1(X))=n$ then there exists a simplicial subdivision
$X'$ of $X$ and a triangulation of $\bigvee_1^n S^1,$ such that
there exists a simplicial map 
$f:X\to \bigvee_1^n S^1$ inducing an epimorphism $\pi_1(X)\to 
\pi_1(\bigvee_1^n S^1).$ Let $\beta_i\in C^1(\bigvee_1^n S^1)$
be a cocycle assigning $1$ to a certain $1$-simplex lying in the
$i$th circle and $0$ to all other simplices of $\bigvee_1^n S^1.$
Then $\alpha_i=f^*(\beta_i)\in C^1(X'),$ for $i=1,...,n,$ satisfy the
conditions of the definition of $c(X)$ and, hence,
$c(X)\geq d(\pi_1(X)).$

Conversely, let $\alpha_1,...,\alpha_n\in C^1(X')$ satisfy the
conditions of the definition of $c(X).$ Consider the map
$f:(X')^1\to \bigvee_1^n S^1$ which maps all vertices of $X'$ to the 
base point of $\bigvee_1^n S^1$
and any $1$-simplex $v$ to a curve $f(v)$ in $\bigvee_1^n S^1$
going $\alpha_1(v)$ times around the first circle,
$\alpha_2(v)$ times around the second, and so on (in that order).

Since $\alpha_i\cup \alpha_j=0$ in $C^2(X')$ for $i\ne j$
for any $2$-simplex $<x_0,x_1,x_2>$ in $X'$
one of the following $3$ conditions holds:
either $\alpha_i(<x_0,x_1>)=\alpha_i(<x_1,x_2>)=0$ for all 
$i$ different than
a certain $k,$ or $\alpha_i(<x_0,x_1>)=0$ for all $i$
or $\alpha_i(<x_1,x_2>)=0$ for all $i.$ 
In each of these cases $f$ extends onto $<x_0,x_1,x_2>.$
Hence, $f$ extends onto the $2$-skeleton, $(X')^2,$ and consequently
it extends to a map $f:X'\to \bigvee_1^n S^1.$ Since 
$[\alpha_1],...,[\alpha_n]\in 
H^1(X,\Z)$ are linearly independent, the image of the map 
$f_*: H_1(X,\Z)\to H_1(\bigvee_1^n S^1)$ is isomorphic with $\Z^n.$
Therefore  the image of $f_*: \pi_1(X)\to \pi_1(\bigvee_1^n S^1)$ 
is the free group of rank $n.$ Hence $d(\pi_1(X))\geq c(X).$
\end{proof}

We will say that a set of vectors $v_1,...,v_n\in H^1(X,\Z)$
is primitive if the map $(v_1,...,v_n): H_1(X,\Z)\to \Z^n$ is an
epimorphism. Primitive vectors are linearly independent.
For any topological space $X,$ let $c_k(X)$ for $k=1,2,...,$
be the cardinality of the largest set of primitive vectors
$v_1,...,v_n\in H^1(X,\Z)$ such that the all Massey products
$<v_{i_1},...,v_{i_l}>$ for $l\leq k$ are uniquely defined and equal to $0.$ 
The definition of the Massey product
can be found for example in \cite{Fenn}.
In particular, $c_2(X)$ is the cardinality of the largest
primitive set of vectors $v_1,...,v_n\in H^1(X,\Z)$ such that 
$v_i\cup v_j=0$ for
all $i,j.$ We will call $c_2(X)$ the algebraic cut number of $X.$

If $\alpha_1,...,\alpha_n\in C^1(X)$ are as in Definition \ref{simplical_cut}
then all Massey products $<[\alpha_{i_1}],...,[\alpha_{i_k}]>$ are $0.$
For that reason we introduce the notation $c_{\infty}(X)=c(X).$ 
We will see in Corollary \ref{c_infty<c_k} that 
$c_{\infty}(X)\leq c_k(X)$ for any $k$
and therefore
\begin{equation}
\label{c-ineq}
c_\infty(X)\leq ... \leq c_k(X)\leq ... \leq c_2(X)\leq c_1(X)=b_1(X).
\end{equation}
Notice that the numbers $c_k(X)$ depend on the homotopy type of $X$ only.
\footnote{One uses here \cite[Lemma 6.2.5]{Fenn}.}
Furthermore, the next proposition shows that these numbers can be derived
from $\pi_1(X).$ Let $K(\Gamma,1)$ denote the Eilenberg-MacLane space for 
$\Gamma$ and let $c_k(\Gamma)=c_k(K(\Gamma,1)).$

\begin{proposition}
If $X$ is a CW-complex then
$c_k(X)=c_k(\pi_1(X))$ for all $k.$
\footnote{Because of the assumption of Definition \ref{simplical_cut},
for $k=\infty$ we assume that $X$ is 
a simplicial complex.}
\end{proposition}

\begin{proof} By Proposition \ref{c=d}, for $k=\infty$ we have 
$c_\infty(X)=d(\pi_1(X))= c_\infty(K(\pi_1(X),1)).$  Hence, we may assume that 
$k<\infty.$ The space $K(\pi_1(X),1)$ can be constructed 
from $X$ by attaching $n$-cells, $n\geq 3.$ Therefore the natural map
$\psi: X\to K(\pi_1(X),1)$ yields a one-to-one correspondence between 
the $i$-cells of $X$ and the $i$-cells of $K(\pi_1(X),1)$ for $i=0,1,2$ 
and it induces an isomorphism $H^1(K(\pi_1(X),1),\Z)\to H^1(X,\Z)$ and 
a monomorphism 
$H^2(K(\pi_1(X),1),\Z)\to H^2(X,\Z).$ Now the statement follows directly 
from the definition of the Massey product (defined for cell cohomology).
\end{proof}

\begin{proposition}\label{ineq_for_epi}
If $G_1$ is an epimorphic image of $G_0$ then
$c_k(G_1)\leq c_k(G_0)$ for each $k.$
\end{proposition}

\begin{proof} For  $k=\infty$ the statement follows from Proposition
\ref{c=d} and therefore assume that $k<\infty.$ 
Since $G_1$ is an epimorphic image of $G_0,$ the Eilenberg-MacLane spaces for
$G_0$ and $G_1$ may be constructed in such a way that there exists 
a cellular map $f:K(G_0,1)\to K(G_1,1)$ which is a bijection
between the $1$-cells. (Consequently, $f_*: G_0=\pi_1(K(G_0,1))\to 
\pi_1(K(G_1,1))=G_1$ is an epimorphism). Let $c_k(G_1)=n$ and let
$\{v_1,...,v_n\}$ be a primitive set of elements of $H^1(G_1,\Z)$ such that
all Massey products $<v_{i_1},...,v_{i_k}>$ are $0.$
Then $\{f^*(v_1),...,f^*(v_n)\}$ is a primitive set of elements of 
$H^1(G_0,\Z)$ and we claim that all Massey products 
$<f^*(v_{i_1}),...,f^*(v_{i_k})>$ are $0.$ Indeed, since
$f$ is a bijection between $C_1(K(G_0,1))$ and $C_1(K(G_1,1))$
any defining set for the Massey product $<f^*(v_{i_1}),...,f^*(v_{i_k})>$
is of the form $(f^*(a_{ij}))_{i,j}$ where $(a_{ij})_{i,j}$ is
a defining set for $\{v_1,...,v_n\}.$ Hence 
$<f^*(v_{i_1}),...,f^*(v_{i_k})>=\mbox{$f^*(<v_{i_1},...,v_{i_k}>)$}=0.$
\end{proof}

By Proposition \ref{c=d}, $c_\infty(\Gamma)=d(\Gamma).$ Therefore,
it is natural to ask whether there exists a similar
interpretation of the numbers $c_k(\Gamma)$ for $k<\infty.$
One may consider a sequence of numbers $d_k(\Gamma)$ defined as follows:
If $F_{n,k}$ denotes the $k$th group in the lower central series of $F_n,$
$F_{n,0}=F_n, F_{n,k+1}=[F_{n,k},F_n],$ then
$d_k(\Gamma)$ is the largest $n$ for which there exists an
epimorphism from $\Gamma$ onto $F_n/F_{n,k}.$ Let $d_\infty(\Gamma)=d(\Gamma).$
The following sequence 
of inequalities parallels (\ref{c-ineq}):

$$d_\infty(\Gamma)\leq ...\leq d_k(\Gamma)\leq ... \leq d_1(\Gamma).$$

\begin{conjecture} \label{conj} $c_k(\Gamma)= d_k(\Gamma)$ for any $k.$
\footnote{One needs to show that if $c_k(\Gamma)=n$ then there exists
an epimorphism $\Gamma\to F_n/F_{n,k}.$
Here is one possible approach suggested to us: The Magnus expansion
of the free group $F_n$ provides an embedding of $F_n$ into the
ring of formal power series in $n$ non-commuting variables, \cite{MKS}.
By taking a quotient of this ring by the ideal generated by all 
monomials of degree $k,$ this map factors to an embedding of $F_n/F_{n,k}$
into the ring of polynomials $\sum_{r=0}^k \sum_{i_1,...,i_r} 
a_{i_1...i_r}x_{i_1}...x_{i_r},$ where $a_{i_1...i_r}\in \Z$ and
$x_1,..., x_n,$ are non-commuting variables.
Therefore, any epimorphism $\Gamma\to F_n/F_{n,k}$ is determined by
the functions $a_{i_1...i_r}:\Gamma\to \Z.$
Those functions are related to the Massey products for $\Gamma.$
In the case of complements of links in homology spheres, this connection 
is mentioned in \cite{T1}.}
\end{conjecture}

Indeed, $c_k(\Gamma)=d_k(\Gamma)$ for $k=1$ and $k=\infty.$
Later on, we will prove it for $k=2$ as well.
Furthermore, we have the following inequality:

\begin{proposition}\label{d<c}
$d_k(\Gamma)\leq c_k(\Gamma)$ for any $k.$ 
\end{proposition}

\begin{proof}
Since for $k=\infty$ this follows from Proposition \ref{c=d}, assume that
$k<\infty.$ If $d_k(\Gamma)=n$ then $F_{n,k}$ is an epimorphic image
of $\Gamma$ and by Proposition \ref{ineq_for_epi}, $c_k(\Gamma)\geq 
c_k(F_{n,k}).$ Now the statement follows from the lemma below.
\end{proof}

\begin{lemma}
$c_k(F_n/F_{n,k})=n.$
\end{lemma}

\begin{proof}
Since $H^1(F_n/F_{n,k},\Z)=\Z^n,$ it is enough to prove that
for all $v_1,...,v_k\in H^1(F_n/F_{n,k}),$ $<v_1,...,v_k>=0.$
Let $U(n)\subset GL(\Z,n)$ be the group of upper triangular $n\times n$ 
matrices all of whose diagonal entries are $1.$ 
By \cite[Thm 2.4]{Dw}, it is sufficient to prove that for any 
$v_1,...,v_k: F_n/F_{n,k}\to \Z$ there exists a homomorphism 
$f:F_n/F_{n,k}\to U(k+1)$ whose $(i,i+1)$ component is $v_i,$ for $i=1,...,k.$
Since $F_n$ is a free group, there exists a homomorphism
$\bar f:F_n\to U(k+1)$ with the above property.
Since the $k$th commutator of $U(k+1)$ is the trivial group, $\bar f$ 
factors through $F_{n,k}$ yielding the required homomorphism.
\end{proof}

Since $c_\infty(X)=d_\infty(\pi_1(X))\leq d_k(\pi_1(X))\leq c_k(\pi_1(X))$
we get

\begin{corollary} \label{c_infty<c_k}
For any simplicial complex $X$ and any $k,$
$c_\infty(X)\leq c_k(X).$\footnote{This statement does not follow
directly from the definition of $c_k(X)$ since
the cohomology classes in Definition \ref{simplical_cut} are not
assumed to be primitive.}
\end{corollary}

Now we are going to prove Conjecture \ref{conj} for $k=2.$
First, we recall a few results from homological algebra.

Given a group $G$ and an abelian group $A,$
the extensions $E$ of $G$ by $A,$ given by the following short exact sequence
$$0\to A\to E\stackrel{\pi}{\to} G\to 0,$$
(up to an isomorphism) are in one-to-one correspondence with elements of 
$H^2(G,A),$ where the action of $G$ on $A$ comes from the action of $E$ 
by conjugation on $A.$
We denote the element of $H^2(G,A)$ corresponding to such extension by
$\xi.$ We say that the sequence splits if there exists a homomorphism 
$\tau: G\to E$ such that $\pi\tau=id_G.$ This happens if and only if $\xi=0.$

\begin{lemma}\label{lift}
(1) A group homomorphism $f:\Gamma \to G$ lifts to a homomorphism
$\bar{f}:\Gamma\to E$ if and only if $f^*(\xi)=0$ in $H^2(\Gamma, A).$\\
(2) If $f^*(\xi)=0$ and $f$ is an epimorphism then
also $\bar{f}$ is an epimorphism.
\end{lemma}

\begin{proof}
(1) Let $E'$ be the pullback of $E$ by $f,$
$$E'=\{(e,\gamma)\in E\times \Gamma: \pi(e)=f(\gamma)\}.$$
The following diagram
\begin{equation}\label{ext}
\begin{array}{ccccccccc}
0 & \to & A & \to & E' & \to & \Gamma& \to & 0\\
  &   &\| & & \downarrow & & \downarrow f & & \\
0 & \to & A & \to & E & \stackrel{\pi}{\to} & \Gamma& \to & 0\\
\end{array}
\end{equation}
is exact and $f$ lifts to $\bar{f}:\Gamma\to E$ if and only if
the top row in (\ref{ext}) splits. By Exercise 6.6.4 in \cite{Weibel} 
this happens if and only if $f^*(\xi)=0.$
We get the second part of the statement of Lemma  \ref{lift}
by chasing the arrows of diagram (\ref{ext}).
\end{proof}

\begin{lemma}
Let $g_1,...,g_n$ be the free generators of $F_n$
and let $e_1,...,e_n$ be the canonical basis of $\Z^n.$
There exists an isomorphism $\kappa: \bigwedge^2 \Z^n \to F_{n,1}/F_{n,2}$
such that 
\begin{equation}\label{e1}
\kappa(e_i\wedge e_j)=[g_i,g_j].
\end{equation}
\end{lemma}

\begin{proof}
Since $F_{n,1}/F_{n,2}$ is abelian and $\bigwedge^2 \Z^n$ is a free
abelian group with basis $e_i\wedge e_j,$ $i<j$, there is unique
homomorphism $\kappa: \bigwedge^2 \Z^n \to F_{n,1}/F_{n,2}$
defined by (\ref{e1}) for $i<j.$ Now note that (\ref{e1}) holds for all
$i,j$ and that $\kappa$ is an epimorphism.
Since, by Corollary 5.12(iv) and Theorem 5.11 in \cite{MKS}
$F_{n,1}/F_{n,2}$ is a free abelian group of the same rank as 
$\bigwedge^2 \Z^n,$ $\kappa$ is an isomorphism.
\end{proof}

\begin{theorem}\label{d_2=c_2}
$d_2(\Gamma)= c_2(\Gamma)$ for any group $\Gamma.$
\end{theorem}

\begin{proof}
If $c_2(\Gamma)=n$ then there exists an epimorphism
$f=(f_1,...,f_n):\Gamma\to \Z^n$ such that
$f_i\cup f_j=0$ in $H^2(\Gamma,\Z)$ for any $i,j.$
Consider the central extension
$$0\to \bigwedge\nolimits^2 \Z^n=F_{n,1}/F_{n,2}\to F_n/F_{n,2}
\stackrel{\pi}{\to} F_n/F_{n,1}=\Z^n\to 0.$$
We will show that $f$ lifts to an epimorphism 
$\bar{f}:\Gamma\to F_n/F_{n,2}$ such that
$f=\pi\bar{f}.$

Since $F_{n,1}/F_{n,2}$ lies in the center of $F_n/F_{n,2},$
the above extension corresponds to $\xi\in H^2(\Z^n,\bigwedge^2 \Z^n),$
where the action of $\Z^n$ on $\bigwedge^2 \Z^n$ is trivial.
Consider the function $\sigma: \Z^n\to F_n/F_{n,2},$
$\sigma(\sum_{i=1}^n a_n e_n)=\prod_{i=1}^n g_i^{a_i}.$ The cohomology
class $\xi$ is represented by the cocycle
$\xi(v_1,v_2)=\sigma(v_1)\sigma(v_2)\sigma(v_1v_2)^{-1},$  compare
\cite[\S6.6]{Weibel}. We have $\xi(e_i,e_j)=1$ for
$i\leq j$ and $\xi(e_i,e_j)=[g_i,g_j]$ for $i>j.$
Hence, using the identities of \cite[Thm 5.3]{MKS}
and the isomorphism $\kappa: \bigwedge^2 \Z^n\to F_{n,1}/F_{n,2},$
we see that $\xi$ corresponds to the cocycle
$$\xi(\sum a_ie_i,\sum b_ie_i)=\sum_{i>j} a_ib_j e_i\wedge e_j$$
in $H^2(\Z^n,\bigwedge^2 \Z^n).$
Hence, for $\gamma_1,\gamma_2\in \Gamma,$
$$f^*(\xi)(\gamma_1,\gamma_2)=
\sum_{i>j} f_i(\gamma_1)f_j(\gamma_2)e_i\wedge e_j=
-\sum_{i>j} (f_i\cup f_j)(\gamma_1,\gamma_2)e_i\wedge e_j,$$
compare \cite[V\S 3]{Brown}.
Since $f_i\cup f_j=0$ in $H^2(\Gamma,\Z),$ also $\xi=0$ in
$H^2(\Gamma,\bigwedge^2 \Z^n).$ Now, by Lemma \ref{lift}, 
$f$ extends to an epimorphism $\bar{f}: \Gamma\to F_n/F_{n,2}.$ 
\end{proof}

\section{Cut numbers and Lie algebras}\label{s_lie}
In this section we prove Theorem \ref{thm-lie}. The background 
information about semisimple Lie algebras can be found in \cite{Hu,Sa}.
For any Lie algebra $\g$
let $\Psi_{\g}(x,y,z)=\kappa(x,[y,z]),$ for $x,y,z\in \g,$ 
where $\kappa$ is the Killing form on $\g.$ The form $\Psi_{\g}$ is 
skew-symmetric.

For any closed, oriented $3$-manifold $M,$ let $\Psi_M: \bigwedge^3 
H^1(M,\Q)\to \Q,$ $\Psi_M(x,y,z)=(x\cup y\cup z)[M].$
Whenever convenient, we will assume that $\Psi_M$ is a $3$-form defined on 
$H^1(M,\R).$

\begin{proposition}
Let $M$ be a $3$-manifold and let
$\g$ be a semisimple Lie algebra of compact type over $\R.$ 
If there are isomorphisms of vector spaces $H^1(M,\R)\simeq \R^n \simeq \g$
such that $\Psi_M=\Psi_{\g}: \bigwedge^3 \R^n\to \R$ then
$b_1(M)=dim\, \g$ and $c_2(M)\leq rank\, \g.$
\end{proposition}

\begin{proof} The equality $b_1(M)=dim\, \g$ is obvious.
If $c_2(M)=k$ then there exist linearly independent vectors
$v_1,...,v_k\in H^1(M,\R)$ such that\\ $\Psi_M(v_i,v_j,H^1(M,\R))=0,$
for any $i,j.$ 
Hence, under the identification\\ $H^1(M,\R)\simeq \g,$
the vectors $v_1,...,v_k$ span a $k$-dimensional abelian subalgebra of $\g.$
Since $\g$ is of compact type, any two maximal abelian Lie subalgebras
of $\g$ are conjugate and their dimensions are equal to $rank\, g.$ Therefore
$k\leq rank\, \g.$
\end{proof}

Now we are ready to prove Theorem \ref{thm-lie}:
For any complex semisimple Lie algebra $\g$ consider its
Chevalley's basis $X_{\alpha}, H_{\alpha},$ where $\alpha'$s
are positive roots of $\g.$
The Lie algebra $\g_0$ spanned by $iH_\alpha, X_\alpha-X_{-\alpha},
i(X_\alpha+X_{-\alpha})$ is a real compact form of $\g,$ compare 
\cite[Ch. I\S 10]{Sa}.
Furthermore, since in Chevalley's basis $\g$ has integral structural constants,
$\Psi_{\g_0}$ is also defined over $\Z.$ Hence $\Psi_{\g_0}: 
\bigwedge^3 \Q^n\to \Q.$
By Sullivan's theorem (\cite{Su}) there exists a $3$-manifold
$M$ with $H^1(M,\Q)\simeq \Q^n$ and $\Psi_M=\Psi_{\g_0}: \bigwedge^3
\Q^n\to \Q.$\Box\\

For any $\g,$ a manifold $M$ as above can be explicitly constructed using
the method of Sullivan, \cite{Su}.

\section{Intermezzo: classification of skew-symmetric forms}
\label{class_forms}

Let $\K$ be a field of characteristic $\chi(\K)\ne 2.$
For any $\Psi: \bigwedge^s \K^n\to \K^m$ we call the
$k$-nullity of $\Psi,$ $null_k(\Psi),$
the maximum of dimensions of vector subspaces $W\subset \K^n$ such that
$$\Psi(\underbrace{W,...,W}_k, \underbrace{\K^n,...,\K^n}_{s-k})=0.$$
$n-null_1(\Psi)$ is called the rank of $\Psi$ and it is equal to the
minimum of dimensions of spaces $V$ such that $\Psi$ factors through
$\bigwedge^s V.$ We will call $null_2(\Psi)$ the
nullity of $\Psi.$ 

The natural $GL_n(\K)$-action on $\K^n$ yields a $GL_n(\K)$-action on forms\\
$\bigwedge^s \K^n\to \K^m,$ $gf(x_1,...,x_n)=f(g^{-1}x_1,...,g^{-1}x_n)$ for
$g\in GL_n(\K),$ $x_1,...,x_n\in \K.$
We say that $\Psi, \Psi':\bigwedge^s \K^n\to \K^m$ are isomorphic if they
belong to the same orbit of the $GL_n(\K)$-action.

The argument used in the previous section implies the following

\begin{fact}
For any $\Psi: \bigwedge^3 \Q^n\to \Q$ with nullity $c$ there exists a 
$3$-manifold $M$ with $b_1(M)=n$ and with $c_2(M)\leq c.$
\end{fact}

Therefore, one is interested in constructing $3$-forms of small nullity. 
This simple looking problem is not easy to 
approach by any straightforward method. A part of the problem is that
$3$-forms do not admit classification similar to that for
$2$-forms: each skew-symmetric $2$-form
is uniquely determined (up to isomorphism) by its rank. More precisely,
for any $\Phi: \bigwedge^2 \K^n\to \K$ of rank $r$ there exists a basis 
$e_1,...,e_n$ of $\K^n$ such that $\Phi$ in the dual basis is given  by
$$\Phi=e^1\wedge e^{r+1}+ ... + e^r\wedge e^{2r}.$$

In contrast to $2$-forms, $3$-forms are not classified by 
their ranks.
The classification of $3$-forms of rank at most $8$ over algebraically
closed fields can be found in \cite[\S35]{Gu}. There are no $3$-forms 
of rank $1,2,4.$ 
There is only one $3$-form (up to an isomorphism) of rank $3,$
$e^1\wedge e^2\wedge e^3,$ and only one $3$-form of rank
$5,$ $e^1\wedge e^2\wedge e^3+e^1\wedge e^2\wedge e^4.$
There are two types of forms of rank $6:$
$$e^1\wedge e^2\wedge e^3+e^4\wedge e^5\wedge e^6 \text{\ and\ }
e^1\wedge e^2\wedge e^3+e^3\wedge e^4\wedge e^5+
e^2\wedge e^5\wedge e^6.$$
The classification of $3$-forms of rank $7,8$ is more complicated.
The classification of $3$-forms of rank $9$ (over $\C$ only)
was carried out in \cite{VE}. 
The classification of $3$-forms of rank larger than $9$ is not known.

\begin{proposition}\label{non-isom}
For any infinite field $\K$ there are infinitely many non-isomorphic forms
$\Psi:\bigwedge^3 \K^n\to \K$ for $n\geq 9.$
\end{proposition}

\begin{proof}
Assume that $n\geq 9$ and that up to an isomorphism
$\Psi_1,...,\Psi_k$ are the only skew-symmetric $3$-forms on $\K^n.$
Let $X$ be the algebraic 
closure of $\bigcup_i GL_n(\o{\K})\Psi_i$ in 
$(\bigwedge^3 \o{\K})^*.$ 
We have 
$dim\, GL_n(\o{\K})=n^2<{n \choose 3}=dim\, (\bigwedge^3 \o{\K})^*$
and, therefore, $dim\, X<dim\, (\bigwedge^3 \o{\K})^*.$
By Lemma \ref{open}(2) (applied to $U=(\bigwedge^3 \o{\K})^*\setminus X$) 
there exists $\Psi\in (\bigwedge^3 \K)^*$
such that $\Psi\not\in X.$ This contradicts the initial assumption.
\end{proof}

\begin{lemma}\label{open} Let $\K$ be an infinite field.
(1) For any non-zero $f\in \o{\K}[x_1,...,x_n]$ there exist 
infinitely many points $(a_1,...,a_n)\in \K^n$ such that 
$f(a_1,...,a_n)\ne 0.$\\
(2) If $U$ is a (Zariski) open subset of $\o{\K}^n,$
then there are infinitely many points in $U$ defined over 
$\K$ (ie. with all their coordinates in $\K$).
\end{lemma}

\begin{proof} We prove (1) by induction. Assume that the statement holds
for $n-1.$ Any non-zero $f\in \o{\K}[x_1,...,x_n]$ is either
constant (and then the statement of lemma holds) or it has a positive 
degree with respect to at least one variable, say $x_1.$ 
Therefore $f=\sum x_1^i f_i,$ where $f_i\in 
\o{\K}[x_2,...,x_n]$ and $f_k\ne 0$ for at least one $k>0.$
By our assumption there exist $a_2,...,a_n\in \K$ such that
$f(a_2,...,a_n)\ne 0$ and hence, $f(\cdot, a_2,...,a_n)$ is a 
non-zero polynomial in $\o{\K}[x_1].$ Since $f(\cdot, a_2,...,a_n)$
has only finitely many zeros, $f(a_1, a_2,...,a_n)\ne 0$ for infinitely
many $a_1\in \K.$ 

(2) Let $f\in \o{\K}[x_1,...,x_n]$ be one of the defining polynomials of 
$X=\o{\K}^n\setminus U,$ ie. let the zero set of $f$ contain $X.$ 
By (1) there are infinitely many points
$(a_1,...,a_n)\in \K^n$ such that $f(a_1,...,a_n)\ne 0.$
Each of these points lies in $U.$
\end{proof}

The difficulties in understanding $3$-forms explain the reason
for using semisimple Lie algebras to construct $3$-forms of small nullity
in the previous section.

\section{Nullity of skew-symmetric forms}
Since skew-symmetric forms cannot be classified directly, in this section 
we study their properties in more subtle ways. Our goal is proving that if 
$\K$ is an algebraically closed field then (a) for any $n$ there exists 
$\Psi: \bigwedge^3 \K^n\to \K$ with $null(\Psi)\leq
2$ and (b) for any $k$ the set of forms $\Psi: \bigwedge^3 \K^n\to \K$ with
$null(\Psi)\geq k$ forms a closed algebraic subset of the space of
all $3$-forms. In the next section we will use these two results to
prove Theorem \ref{main}.

Since our methods of proof are not specifically limited to
$3$-forms, we formulate our results more broadly in terms of arbitrary
skew-symmetric forms, hoping that one day they will 
find new applications, unrelated to cut numbers of $3$-manifolds.

We start with the following result due to T. Goodwillie.

\begin{theorem}\label{Goodwillie}
Let $\K$ be algebraically closed.
(1) If $m\leq 2n-4$ then the nullity of any form 
$\Phi: \bigwedge^2 \K^n\to \K^m$ is at least $2.$
(2) If $m> 2n-4$ then there exists
$\Phi: \bigwedge^2 \K^n\to \K^m$ with $null(\Phi)=1.$ 
\end{theorem}

\begin{proof}
Consider the Grassmannian of $2$-planes in $\K^n,$ $G_2(\K^n),$ embedded
into $\P(\bigwedge^2 \K^n)$ by the map $V^2\to V^2\wedge V^2.$
Observe that $null(\Phi)\geq 2$ if and only
if $\P(Ker\, \Phi)\cap G_2(\K^n)\ne \emptyset$ in $\P(\bigwedge^2 \K^n).$
A projective subvariety of dimension $s$ of $\P(\K^N)$
intersects every projective subspace of codimension at most $s$ in 
$\P(\K^N),$ but it does not intersect some projective subspaces of codimension 
$s+1;$ compare \cite[Thm I.7.2]{Ha}. Since $dim\, G_2(\K^n)=2n-4$ and
$codim\, P(Ker\, \Phi)=m,$  $null(\Phi)\geq 2$ for all
$\Phi:\bigwedge^2 \K^n\to \K^m$ with $m\leq 2n-4.$

If $m\geq 2n-4$ then there exists a plane $W\subset \bigwedge^2 \K^n$
of codimension $m$ such that $\P(W)\cap G_2(\K^n)=\emptyset.$
For any such $W,$ the $2$-form $\Phi: \bigwedge^2 \K^n\to \bigwedge^2 
\K^n/W=\K^m$ has nullity $1.$
\end{proof}

By considering a $3$-form $\Psi:  \bigwedge^3 \K^n\to \K$ as a $2$-form,
$\bigwedge^2 \K^n\to (\K^n)^*$ we get the following corollary:

\begin{corollary}\label{Good_cor}
For $n\geq 4$ the nullity of every $\Psi: \bigwedge^3 \K^n\to \K$ is at 
least $2.$
\end{corollary}

Theorem \ref{Goodwillie} does not hold over fields which are not 
algebraically closed. Nevertheless, it is very possible that the 
above corollary
holds over any field $\K$ of characteristic $\chi(\K)\ne 2.$ 
\vspace*{.1in}

\noindent {\bf Conjecture:} {\em For any field $\K$ of characteristic 
$\chi(\K)\ne 2$ and for $n\geq 4,$ the nullity of every
$\Psi: \bigwedge^3 \K^n\to \K$ is at least $2.$}

\begin{proposition}\label{even}
The conjecture holds for even $n.$\footnote{We can also prove the 
conjecture for $n=5.$}
\end{proposition}

\begin{proof} For any $v\in \K^n$ there is an $n\times n$ matrix $A(v)$
such that $\Psi(v,w_1,w_2)=w_1^{T}\cdot A\cdot w_2$ for any $w_1,w_2\in \K^n.$
Since $A\cdot v=0,$ $rank\, A\leq n-1.$ Since $\chi(\K)\ne 2$ and
$A$ is skew-symmetric,
the null space of $A,$ $\{w: Aw=0\},$ has an even dimension and, hence, 
there is $w$ linearly independent from $v$ such that $\Psi(v,w,\K^n)=0.$
\end{proof}

\begin{theorem}\label{closed} Let $\K$ be an algebraically closed field.
For any $s,k,r,n,m,$ the set of forms 
$\Psi: \bigwedge^s \K^n\to \K^m$ with $null_k(\Psi)\geq r$
is a closed algebraic subset of $(\bigwedge^s \K^n)^*\otimes \K^m$ 
defined by homogeneous polynomials.\footnote{Recall that
$null_k(\Psi)$ was defined in Section \ref{class_forms}.}
\end{theorem}

Fix $s,k,r,n,m.$
For each $r$-dimensional subspace $V\subset \K^n$ consider
$$\imath(V)=\bigwedge\nolimits^k V\wedge \bigwedge\nolimits^{s-k} \K^n=
Span\, \{v_1\wedge ... \wedge v_k\wedge w_1\wedge ...\wedge w_{s-k}\}
\subset \bigwedge\nolimits^s \K^n,$$ where $v_1,..., v_k\in V,$ 
$w_1,...,w_{s-k}\in \K^n.$ 

\begin{lemma}\label{dim}
For any $V$ as above, $dim\, \imath(V)=d,$ where
$$d=\sum_{k\leq i\leq r}
{i-1\choose k-1}{n-i\choose s-k}.$$
\end{lemma}

\begin{proof}
If $e_1,...,e_n$ is a basis of $\K^n$ such that $V=Span\{e_1,...,e_r\}$ 
then $\imath(V)$ has a basis composed of elements 
$e_{i_1}\wedge ... \wedge e_{i_k}\wedge ... \wedge e_{i_s},$
where $1\leq i_1< ...< i_k< ...< i_s\leq n$ and
$i_k\leq r.$ The number of basis elements of $\imath(V)$ with
$i_k=i$ is ${i-1\choose k-1}{n-i\choose s-k}.$
\end{proof}

Let $$X_{k,r,s,n}=\{\imath(V): V\in G_r(\K^n)\}\subset 
G_d(\bigwedge\nolimits^s \K^n).$$
We have
\begin{equation}\label{iff}
null_k(\Psi)\geq r \text{\quad if and only if\quad }
G_d(Ker\, \Psi)\cap X_{k,r,s,n}\ne \emptyset \text{\ in\ 
$G_d(\bigwedge\nolimits^s \K^n).$}
\end{equation}

\noindent{\bf Proof of Theorem \ref{closed}:}
The proof relies on elimination theory for projective spaces, see
\cite[Ch. 8\S 5]{CLO}. (Compare also \cite[\S80, p.8]{vW}).

Let $\{e_1,...,e_N\}$ be a basis of
$\bigwedge\nolimits^s \K^n,$ $N={n \choose s}.$
Using this basis we identify
$\bigwedge\nolimits^s \K^n$ with $\K^N.$
The Grassmannian  $G_d(\K^N)$ embeds into a projective space by 
$P:G_d(\K^N)\hookrightarrow \P(\bigwedge\nolimits^d \K^N),$
$P(W)=\bigwedge\nolimits^d W.$
Each $P(W)$ can be written as $$\sum_{1\leq i_1<...<i_d\leq N}
P_{i_1...i_d}(W)e_{i_1}\wedge ... \wedge e_{i_d},$$
where the numbers $P_{i_1...i_d}(W),$ called the Pl\"ucker 
coordinates of $W,$ are defined up to common scalar
multiplication.

Let $\Psi:\K^N \to \K^m$ be given by
as $\Psi(x)=\sum_{i=1}^N \Psi_i e^i(x),$
where $\{e^1,...,e^N\}$ is the basis of $(\K^N)^*$ dual to $\{e_1,...,e_N\}$
and $\Psi_i=(\Psi_{ij})_{1\leq j\leq m}\in \K^m.$

Let $M_{kl}(\K)$ denote the set of $k\times l$ matrices with entries in $\K.$

\begin{lemma}\label{elimination}
There exist polynomial functions $F_1,...,F_f: \bigwedge\nolimits^d \K^N
\times M_{N,m}(\K) \to \K$ such that
for any $\Psi=\sum_{i=1}^N \Psi_ie^i,$
the vector $\sum_{1\leq i_1<...< i_d\leq N}
P_{i_1...i_d}e_{i_1}\wedge ...\wedge e_{i_d}
\in \bigwedge\nolimits^d \K^N$
belongs to $\bigwedge\nolimits^d Ker\, \Psi$
if and only if
$$F_i((P_{i_1...i_d})_{1\leq i_1<...< i_d\leq N},
\Psi_1,...,\Psi_N)=0 \text{\ for $1\leq i\leq f$}.$$
Furthermore, the polynomial functions $F_1,...,F_f$ are
homogeneous with respect to the coordinates of $\bigwedge\nolimits^d \K^N.$
\end{lemma}

\begin{proof} Denote $dim\, Ker\, \Psi=N-m$ by $t.$ Since the statement is
obvious for $t<d,$ assume that $t\geq d.$
The vectors $v_1=\sum_{j=1}^N a_{1j}e_j, ..., v_t=\sum_{j=1}^N a_{tj}e_j$ 
belong to $Ker\, \Psi$ if and only if
\begin{equation}\label{eq1}
\sum_{j=1}^N a_{ij}\Psi_j=0 \text{\ in $\K^m$ for every $1\leq i\leq t.$}
\end{equation}
Note that 
$$v_{i_1}\wedge ...\wedge v_{i_d}=\sum_{1\leq j_1< ...<j_d\leq N}
det\left(\begin{matrix}
a_{i_1,j_1} & ... & a_{i_d,j_1}\\
... & ... & ...\\
a_{i_1,j_d} & ... & a_{i_d,j_d}
\end{matrix}\right) e_{j_1}\wedge ...\wedge e_{j_d}.
$$
Therefore, $w\in \bigwedge\nolimits^d \K^N$
belongs to $\bigwedge\nolimits^d Ker\, \Psi$ if and only if there exist
$(a_{ij})\in M_{t,N}$ satisfying (\ref{eq1}) and there exist
$b_{i_1...i_d}\in \K$ for all $1\leq i_1<...<i_d\leq t$
such that 
\begin{equation}\label{eq2}
w=\sum_{1\leq j_1< ...<j_d\leq N\atop 1\leq i_1< ...<i_d\leq t}
b_{i_i...i_d} det\left(\begin{matrix}
a_{i_1,j_1} & ... & a_{i_d,j_1}\\
... & ... & ...\\
a_{i_1,j_d} & ... & a_{i_d,j_d}
\end{matrix}\right) e_{j_1}\wedge ...\wedge e_{j_d}.
\end{equation}

Since (\ref{eq1}), (\ref{eq2}) are homogeneous equations with respect to
$(a_{ij})$ and\\ $(b_{i_1...i_d})_{1\leq i_1<...<i_d\leq t},$
by \cite[Ch. 8\S 5 Thm 8]{CLO}
there exist polynomial functions $H_1,...H_h$ on $\bigwedge\nolimits^d \K^N
\times M_{N,m}(\K)$ such that for any $\Psi_1,...,\Psi_N\in \K^m$ 
the equations
$$H_i(\cdot, \Psi_1,...,\Psi_N)=0 \text{\ for $1\leq i\leq h$}$$ 
describe the set $\bigwedge\nolimits^d Ker\, \Psi.$
For each $i$ let $H_{ij}$ be the sum of monomials of $H_i$ of
total degree $j$ with respect to the coordinates of 
$\bigwedge\nolimits^d \K^N,$ $H_i=\sum_j H_{ij}.$
Since $\alpha\in \bigwedge\nolimits^d Ker\, \Psi$
if and only if $c\alpha\in \bigwedge\nolimits^d Ker\, \Psi,$ for any
$c\in \K^*,$ $H_{ij}(\alpha,\Psi_1,...,\Psi_N)=0$ for any 
$\Psi_1,...,\Psi_N\in \K^m$ and
$\alpha\in \bigwedge\nolimits^d Ker\, \Psi,$
for $\Psi=\sum \Psi_i e^i.$
Hence the polynomials $H_{ij}$ satisfy the conditions of the lemma.
\end{proof}

We continue the proof of Theorem \ref{closed}.
As before, let $\bigwedge^s \K^n=\K^N$ and consider the embedding 
$P:G_d(\K^N)\hookrightarrow \P(\bigwedge\nolimits^d \K^N).$ 
Let $G_1,...,G_g: \bigwedge\nolimits^d \K^N\to \K$ be homogeneous 
polynomials describing $P(X_{k,r,s,n})\subset \P(\bigwedge\nolimits^d \K^N).$
By (\ref{iff}) $null_k(\Psi)\geq r$ if and only if the polynomials
$G_i(\cdot)$ and $F_i(\cdot, \Psi_1,...,\Psi_N)$ (defined in Lemma
\ref{elimination}) have a common solution in
$ \bigwedge\nolimits^d \K^N$ different from $0.$
By \cite[Ch. 8\S 5 Thm 8]{CLO}, the set of $(\Psi_1,...,\Psi_N)$ 
for which such a solution exists is an algebraic set in $M_{N,m}(\K).$
An argument analogous to that given at the end of proof of previous lemma
shows that this set can be described by homogeneous polynomials.
This completes the proof of Theorem \ref{closed}.
\Box

\begin{corollary} If $\K$ is an algebraically closed field then
for any $s,k,r,n,m,$ there is an algorithm (using Gr\"obner bases)
for calculating $null_k(\Psi)$ for $\Psi:\bigwedge^s \K^n\to \K^m.$
\end{corollary}

Any $\Psi: \bigwedge^s \K^n\to \K^m$ yields
$\o{\Psi}: \bigwedge^s \o{\K}^n\to \o{\K}^m$ with
$null_k(\o{\Psi})\geq null_k(\Psi)$ for any $k.$

For any field $\K$ the process of eliminating variables of polynomials 
over $\K$ produces polynomials over the same field $\K.$ Therefore,
the proof of Theorem \ref{closed} implies a more general statement:

\begin{corollary}\label{closed_any_field}
For an arbitrary field $\K$ and for any $s,k,r,n,m,$ the set of forms 
$\Psi: \bigwedge^s \K^n\to \K^m$ with $null_k(\overline{\Psi})\geq r$
is a closed algebraic subset of $(\bigwedge^s \K^n)^*\otimes \K^m$ 
defined by homogeneous polynomials with coefficients in $\K.$
\end{corollary}

\begin{theorem}\label{cut2}
For any algebraically closed field $\K$ and for 
any $n$ there exists $\Psi: \bigwedge^3 \K^n \to \K$ with 
$null_2(\Psi)\leq 2.$
\end{theorem}

\begin{proof}
Since for $n\leq 3$ the statement is obvious, assume that $n\geq 4.$ 
If there exists a hyperplane $W\subset \bigwedge^3 \K^n$ such that 
$G_d(W)\cap X_{k,r,s,n}=\emptyset$ for $k=2,$ $r=s=3,$ then by (\ref{iff}),
the projection map $\Psi: \bigwedge^3 \K^n\to \bigwedge^3 \K^n/W=\K$ will be
a skew-symmetric $3$-form on $\K^n$ with $null_2(\Psi)\leq 2.$
We are going to prove that such $W$ exists.

Let $W_0$ be any hyperplane in $\bigwedge^3 \K^n.$ Since the Grassmannian
of $d$ dimensional planes in $\K^N$ has dimension $dN-d^2$ and
$dim\, X\leq dim\, G_3(\K^n),$ we have
$$dim\, G_d(W_0)+dim\, X \leq
d({n \choose 3}-1)-d^2+ 3n-9$$ and
$$dim\, G_d(\bigwedge\nolimits^3 \K^n)=d{n \choose 3}-d^2.$$
By Lemma \ref{dim}, $d=3n-8$ and, therefore, 
$$dim\, G_d(W_0)+dim\, X< dim\, G_d(\bigwedge\nolimits^3 \K^n).$$
Since $G_d(\bigwedge^3 \K^n)$ is a homogeneous 
$GL(\bigwedge^3 \K^n)$-space,
there exists $g\in GL(\bigwedge^3 \K^n)$ such that $X$ and 
$gG_d(W_0)$ are transversal ie., $X\cap gG_d(W_0)=\emptyset,$ compare 
\cite[Thm III.10.8]{Ha}. (Note that since $G_d(W_0)$ is smooth,
Hartshorne's assumption that $\K$ has characteristic $0$ is not necessary 
here). Since $gG_d(W_0)=G_d(gW_0),$ the hyperplane $W=gW_0$ satisfies the
required condition.
\end{proof}

Note that we proved more:

\begin{corollary} If $\K$ is algebraically closed and $s,k,r,n,m$ 
are such that $r\cdot n-r^2<d\cdot m,$ where $d$ is as in Lemma \ref{dim}, 
then there exists a form $\Psi:\bigwedge\nolimits^s \K^n\to \K^m$
with $null_k(\Psi)< r.$
\end{corollary}

\section{Proof of Theorem \ref{main}} 
\noindent (1) By Corollary \ref{closed_any_field},
$X=\{\Psi: null_2(\o{\Psi})>2\}$ is a closed algebraic subset of
$(\bigwedge^3 \Q^n)^*.$ By Theorem \ref{cut2}, 
$U=(\bigwedge^3 \Q^n)^*\setminus X$ is a nonempty, and hence dense, open
subset of $(\bigwedge^3 \Q^n)^*.$ Since $c(M)\leq null_2(\Psi_M)\leq
null_2(\o{\Psi_M}),$ the set of forms $\Psi_M$ for $3$-manifolds $M$ with 
$c(M)\leq 2$ contains $U.$

\noindent (2) follows from (1).\\
The proof (3) is analogous to the proof of Proposition \ref{non-isom}:
Fix $n\geq 9$ and assume that the $3$-manifolds $M$ with 
$b_1(M)=n$ and $c(M)\leq 2$
represent only finitely many different rational cohomology rings corresponding
to $\Psi_1,...,\Psi_k: \bigwedge^3 \Q^n\to \Q.$
Let $Y=\o{\cup_i GL_n(\C)\Psi_i}\subset (\bigwedge^3 \C^n)^*.$
Since the set $X$ (defined above) and the set $Y$ are closed and strictly 
contained in $(\bigwedge^3 \C^n)^*,$ 
by Lemma \ref{open}(2) there exists $\Psi: \bigwedge^3 \Q^n\to \Q$ such that
$\Psi\not\in X\cup Y.$ A contradiction.
\Box

\end{document}